\documentclass[12pt]{article}
\usepackage{amscd, amssymb, amsmath, calc, enumerate}
\begin{document} 
\renewcommand{\thesubsection}{\arabic{subsection}}
\newenvironment{eq}{\begin{equation}}{\end{equation}}
\newenvironment{proof}{{\bf Proof}:}{\vskip 5mm }
\newenvironment{rem}{{\bf Remark}:}{\vskip 5mm }
\newenvironment{remarks}{{\bf Remarks}:\begin{enumerate}}{\end{enumerate}}
\newenvironment{examples}{{\bf Examples}:\begin{enumerate}}{\end{enumerate}}  
\newtheorem{proposition}{Proposition}[subsection]
\newtheorem{lemma}[proposition]{Lemma}
\newtheorem{definition}[proposition]{Definition}
\newtheorem{theorem}[proposition]{Theorem}
\newtheorem{cor}[proposition]{Corollary}
\newtheorem{conjecture}{Conjecture}
\newtheorem{pretheorem}[proposition]{Pretheorem}
\newtheorem{hypothesis}[proposition]{Hypothesis}
\newtheorem{example}[proposition]{Example}
\newtheorem{remark}[proposition]{Remark}
\newtheorem{ex}[proposition]{Exercise}
\newtheorem{cond}[proposition]{Conditions}
\newtheorem{cons}[proposition]{Construction}
\newcommand{\llabel}[1]{\label{#1}}
\newcommand{\comment}[1]{}
\newcommand{\sr}{\rightarrow}
\newcommand{\lr}{\longrightarrow}
\newcommand{\xr}{\xrightarrow}
\newcommand{\dw}{\downarrow}
\newcommand{\bdl}{\bar{\Delta}}
\newcommand{\zz}{{\bf Z\rm}}
\newcommand{\zq}{{\bf Z}_{qfh}}
\newcommand{\nn}{{\bf N\rm}}
\newcommand{\qq}{{\bf Q\rm}}
\newcommand{\nq}{{\bf N}_{qfh}}
\newcommand{\oo}{\otimes}
\newcommand{\uu}{\underline}
\newcommand{\ih}{\uu{Hom}}
\newcommand{\af}{{\bf A}^1}
\newcommand{\wt}{\widetilde}
\newcommand{\gm}{{\bf G}_m}
\newcommand{\dsr}{\stackrel{\sr}{\scriptstyle\sr}}
\newcommand{\PP}{$P_{\infty}$}
\newcommand{\tp}{\tilde{D}}
\newcommand{\HH}{$H_{\infty}$}
\newcommand{\ii}{\stackrel{\scriptstyle\sim}{\sr}}
\newcommand{\BB}{_{\bullet}}
\newcommand{\D}{\Delta}
\newcommand{\colim}{{\rm co}\hspace{-1mm}\lim}
\newcommand{\cf}{{\it cf} }
\newcommand{\msf}{\mathsf }
\newcommand{\mcal}{\mathcal }
\newcommand{\ep}{\epsilon}
\newcommand{\tl}{\widetilde}
\newcommand{\ub}{\mbox{\rotatebox{90}{$\in$}}}
\newcommand{\piece}{\vskip 3mm\noindent\refstepcounter{proposition}{\bf
\theproposition}\hspace{2mm}}
\newcommand{\subpiece}{\vskip 3mm\noindent\refstepcounter{equation}{\bf\theequation}\hspace{2mm}}{\vskip
3mm}
\numberwithin{equation}{subsection}
\begin{center}
{\Large\bf On the zero slice of the sphere spectrum}\\
\vskip 4mm
{\large Vladimir Voevodsky}\\
{\em December 2002}
\end{center}
\vskip 4mm
\tableofcontents
\subsection{Introduction}
In \cite{Open} we introduced the slice filtration on the motivic
stable homotopy category which is a motivic analog of the filtration
by the subcategories of n-connected spectra in topology. Since the
inclusion functors between different terms of the filtration have
right ajoints it makes sense to speak of the projection $s_n(E)$ of a
spectrum $E$ to the n-th quotient of this filtration. This projection
which is again an object of the motivic stable homotopy category is
called the n-th slice of $E$. Its topological analog is the spectrum
$\Sigma^n H_{\pi_n(E)}$ where $\pi_n(E)$ is the n-th stable homotopy
group of $E$ and $H_A$ is the Eilenberg-MacLane spectrum corresponding
to the abelian group $A$. In this sense slices provide a motivic
replacement of the stable homotopy groups. 

The goal of this paper is to prove over fields of characteristic zero
the main conjecture of \cite{Open} which asserts that the zero slice
of the sphere spectrum ${\bf 1}$ is the motivic Eilenberg-MacLane
spectrum $H_{\zz}$. Using the analogy between slices and stable
homotopy groups one may interpret this result as a motivic version of
the statement that $\pi_0^s(S^0)=\zz$. As an immediate corollary one
gets that the slices of any spectrum are modules over $H_{\zz}$.

We obtain our main result from an unstable statement about the motivic
Eilenber-MacLane spaces $K_n$. We introduce the notion of a
(homotopically) n-thick space and show that on the one hand the
suspension spectrum of any n-thick space belongs to the $n$-th stage
of the slice filtration and on the other hand that for $n>0$ the cone
of the natural map $T^n\sr K_n$ is $(n+1)$-thick.

The restriction on the characteristic of the base field $k$ appears in
our approach twice. Firstly, we use the model of $K_n$ based on the
symmetric powers of $T^n={\bf A}^n/({\bf A}^n-\{0\})$ which is only
valid over fields of characteristic zero. Secondly, in order to
analyze the structure of the d-th symmetric power we need
invertibility of $d!$ in $k$. We expect that the main result of the
paper remains valid in any characteristic even though the intermediate
steps we use do not. 

The first draft of this paper (with some mistakes) was announced on
the Great Lakes K-theory meeting in Evanston in 2001 and I would like
to thank all the listeners for pointing out gaps in the original
argument.

\subsection{Thom spaces of linear representations and their quotients}
We fix a base field $k$ and let $Sch/k$ denote the category of
quasi-projective schemes over $k$. For a finite group scheme $G$ over
$k$ we denote by $G-Sch/k$ the category of $G$-objects in $Sch/k$. In
\cite{delnotes} we defined an analog of the Nisnevich topology on
$G-Sch/k$. We let $(G-Sch/k)_{Nis}$ denote the corresponding site and
$Spc_{\bullet}(G)$ denote the category of pointed simplicial sheaves
on $(G-Sch/k)_{Nis}$. Following
\cite{delnotes} one introduces the class of $\af$-weak equivalences on
$Spc_{\bullet}(G)$ and defines the corresponding $\af$-homotopy
category as the localization with respect to this class. 

We let $\Sigma^1_s$ and $\Sigma^1_T$ denote the suspensions by the
simplicial circle and the sphere $T=\af/(\af-\{0\})$ respectively, and
$\tilde{\Sigma_s}$ denote the unreduced s-suspension:
$$\tilde{\Sigma_s}F=cone(F_+\sr S^0)$$

For a $G$-equivariant vector bundle $E$ over $X$ we let $Th(E)$ denote
the Thom space $E/E-s_0(X)$ where $s_0$ is the zero section. For
$X=Spec(k)$, equivariant vector bundles are linear representations $V$
of $G$ and $Th(V)$ is the object $V/V-\{0\}$. Note that $Th(V\oplus
W)\cong Th(V)\wedge Th(W)$ and $Th(0)=S^0$. Note also that since $V$
is $\af$-contractible as a $G$-scheme there is a natural
$\af$-equivalence of the form
$$Th(V)\cong \tilde{\Sigma_s}(V-\{0\})$$
We let $Quot_G:Spc_{\bullet}(G)\sr Spc_{\bullet}$ denote
the functor which commutes with colimits and such that
$Quot(X_+)=(X/G)_+$ (see \cite[Sec. 5.1]{delnotes}).

In this section we prove several general results about the structure
of the quotients $Quot_G(Th(V))$. For simplicity, we assume in this
section that $G$ is a {\em finite group}. The following lemma is
straightforward.
\begin{lemma}
\llabel{w1}
Let $V_1, V_2$ be linear representations of $G_1, G_2$. Then there is
a natural isomorphism
$$Quot_{G_1\times G_2}(Th(V_1\times V_2))=Quot_{G_1}(V_1)\wedge
Quot_{G_2}(V_2)$$
In particular if ${\cal O}^n$ is the trivial representation of $G$ of
dimension $n$ then 
$$Quot_G(Th(V\oplus {\cal O}^n))=\Sigma_T^nQuot_G(Th(V))$$
\end{lemma}
For a subgroup $H$ in $G$ we let $V^{\ge H}$ denote the closed subset
of $H$-invariant elements in $V$ and by $V^{>H}$ the subset of
elements whose stabilizer is strictly greater than $H$. Note that for
$H_1\ne H_2$ one has $V^{\ge H_1}\cap V^{\ge H_2}\subset V^{>H_1}\cap
V^{>H_2}$ and in particular
$$(V^{\ge H_1}-V^{>H_1})\cap (V^{\ge H_2}-V^{>H_2})=\emptyset$$
For a closed subset $Z$ of $V$ we let $GZ$ denote the orbit of $Z$
i.e. the closed subset $Im(Z\times G\sr V\times G\sr V)$ where the
first map is the closed embedding and the second one is the action. 
\begin{lemma}
\llabel{l1}
Let $H$ be a closed subgroup in $G$ and $N(H)$ the normalizer of
$H$. Then one has  
\begin{eq}\llabel{f1}
\begin{array}{c}
Quot_{G}((V-GV^{>H})/(V- GV^{\ge
H}))=\\\\
=Quot_{N(H)}((V-V^{>H})/(V-V^{\ge H}))
\end{array}
\end{eq}
\end{lemma}
\begin{proof}
Let $A=G/N(H)$ be the set of subgroups adjoint to $H$. The scheme
$GV^{\ge H}-GV^{>H}$ is (non canonically) isomorphic to $(V^{\ge
H}-V^{>H})\times A$ and there is a $G$-equivariant map $GV^{\ge
H}-GV^{>H}\sr A$. Consider the section of the projection $(GV^{\ge
H}-GV^{>H})\times A\sr (GV^{\ge H}-GV^{>H})$ defined by this map.
Since this projection is etale the image of this section is open and
we may consider the closed complement $C$ to this image. Consider the
square:
$$
\begin{CD}
(V-GV^{\ge H})\times A @>>> (V-GV^{>H})\times A - C\\
@VVV @VVpV\\
V-GV^{\ge H} @>j>> V-GV^{>H}
\end{CD}
$$
where the vertical morphisms come from the obvious projections and the
horizontal morphisms from the obvious embeddings. Let us show that
this is an upper distinguished square in $G-Sch/k$. It is clearly a
pull-back square, $p$ is etale and $j$ is an open embedding. It
remains to check that $p^{-1}(GV^{\ge H}-GV^{>H})\sr GV^{\ge
H}-GV^{>H}$ is an isomorphism which follows from our choice of $C$. We
conclude that
$$
(V-GV^{>H})/(V-GV^{\ge H})=(V-GV^{\ge H})\times A/((V-GV^{>H})\times A
- C)
$$
and therefore
$$Quot_{G}((V-GV^{>H})/(V-GV^{\ge H}))=$$
$$=Quot_{G}((V-GV^{\ge H})\times A/((V-GV^{>H})\times A - C))$$
We further have: 
$$Quot_{G}((V-GV^{\ge H})\times A)=Quot_{N(H)}(V-GV^{\ge H})$$
and 
$$Quot_{G}((V-GV^{>H})\times A - C)=Quot_{N(H)}(V-GV^{>H}-(GV^{\ge
H}-V^{\ge H}))$$
Since $GV^{\ge H}-GV^{>H}=(V^{\ge H}-V^{>H})\coprod (GV^{\ge
H}-V^{\ge H}-GV^{>H})$ we have 
$$(V-GV^{>H}-(GV^{\ge
H}-V^{\ge H}))/(V-GV^{\ge H})=(V-V^{>H})/(V-V^{\ge H})$$
and (\ref{f1}) follows.
\end{proof}
\begin{remark}\rm
An analog of Lemma \ref{l1} holds for any finite etale group scheme
$G$ and any $X$ in $G-Sch/k$.
\end{remark}
Let $H$ be a normal subgroup of $G$. Then we may consider the relative
analog of the $Quot$ functor
$$Quot_{G,H}:Spc_{\bullet}(G)\sr Spc_{\bullet}(G/H)$$
which commutes with colimits and such that $Quot(X_+)=(X/H)_+$ where
$X/H$ is considered with the natural action of $G/H$. One verifies
easily that one has $Quot_G=Quot_{G/H}Quot_{G,H}$.
\begin{lemma}
\llabel{l2}
If $H$ is a normal subgroup of $G$ of order prime to $char(k)$ then
one has:
$$Quot_G((V-V^{>H})/(V-V^{\ge H}))=$$
$$=Quot_{G/H}((V^{\ge H}-V^{>H})_+\wedge
Quot_{G,H}(Th(V/V^{\ge H})))$$
\end{lemma}
\begin{proof}
Since 
$$Quot_G=Quot_{G/H}Quot_{G,H},$$
it is sufficient to show that
$Quot_{G,H}((V-V^{>H})/(V-V^{\ge H}))$ is isomorphic to $(V^{\ge
H}-V^{>H})_+\wedge Quot_{G,H}(Th(V/V^{\ge H}))$ as a
$G/H$-space. Since the order of $H$ is prime to $char(k)$ there is an
isomorphism $V=V^{\ge H}\oplus (V/V^{\ge H})$. Using this isomorphism
we get an isomorphism
$$(V-V^{>H})/(V-V^{\ge H})=(V^{\ge H}-V^{>H})_+\wedge Th(V/V^{\ge H})$$
Since the action of $H$ on $V^{\ge H}$ is trivial we get 
$$Quot_{G,H}((V-V^{>H})/(V-V^{\ge H}))=(V^{\ge H}-V^{>H})_+\wedge
Quot_{G,H}(Th(V/V^{\ge H}))$$
\end{proof}
Combining Lemmas \ref{l1} and \ref{l2} we get the following result.
\begin{proposition}
\llabel{basic}
Let $G$ be a finite group of order prime to $char(k)$. Let further $V$
be a linear representation of $G$ over $k$ such that $V^{\ge
G}=0$. Then $Quot_G(Th(V))$ belongs to the smallest class which is
closed under cones, finite coproducts and $\af$-equivalences and 
contains the following objects:
\begin{enumerate}
\item $\tilde{\Sigma_s}(Quot_G(V_0))$ where $V_0$ is the open
subscheme in $V$ where $G$ acts freely
\item $Quot_{N(H)/H}((V^{\ge H}-V^{>H})_+\wedge Quot_{N(H),H}(Th(V/V^{\ge
H})))$ for all subgroups $H$ in $G$ such that $H\ne e,G$.
\end{enumerate}
\end{proposition}
\begin{proof}
Consider the sequence of open embeddings $V_0\sr V_1\sr\dots\sr
V_{N}=V$ where $V_i$ is the subscheme of points with the stabilizer of
order no more than $i+1$ and $N=|G|-1$. In particular $V_0$ is the
open subscheme where $G$ acts freely. Under our assumption 
$V_{N-1}=V-\{0\}$ and hence $Th(V)=\tilde{\Sigma_s}(V_{N-1})$.
\begin{lemma}
\llabel{l4}
Let $X_0\sr X_1\sr\dots\sr X_m$ be a sequence of monomorphisms in
$Spc_{\bullet}(G)$. Then $\tilde{\Sigma_s}(X_m)$ belongs to the
smallest class which is closed under simplicial weak equivalences and
cones and contains $\tilde{\Sigma_s}(X_0)$ and $X_i/X_{i-1}$ for
$i=1,\dots,m$.
\end{lemma}
\begin{proof}
Denote the smallest class satisfying the conditions of the lemma by
$B$. let us show by induction on $m$ that $\tilde{\Sigma_s}(X_m)$ is
in $B$. For $m=0$ the statement is obvious. By induction we may assume
that $\tilde{\Sigma_s}X_{m-1}$ is in $B$. The digram
$$
\begin{CD}
pt @>>> S^0 @>>> S^0\\
@VVV @VVV @VVV\\
X_m/X_{m-1} @>>> \Sigma_s X_{m-1,+} @>>> \Sigma_s X_{m,+}\\
@VVV @VVV @VVV\\
X_m/X_{m-1} @>>> \tilde{\Sigma_s} X_{m-1} @>>> \tilde{\Sigma_s}
X_{m}
\end{CD}
$$
whose rows and columns are cofibration sequences, shows that
$\tilde{\Sigma_s} X_m = cone(X_m/X_{m-1}\sr
\tilde{\Sigma_s} X_{m-1})$ and therefore $\tilde{\Sigma_s} X_m$ is in $B$. 
\end{proof}
In view of Lemma \ref{l4} it is sufficient to show that the quotients
$V_i/V_{i-1}$ for $i=1,\dots,N-1$ belong to the class we
consider. This follows directly from Lemmas \ref{l1} and \ref{l2}.
\end{proof}

\subsection{Thick spaces}
\begin{definition}
\llabel{asdiag}
The class of n-thick objects is the smallest class $A_n$ such that:
\begin{enumerate}
\item $A_n$ is closed under $\af$-equivalences
\item $A_n$ is closed under filtering colimits
\item if $(F_{i})$ is a simplicial object in
$\Delta^{op}Spc_{\bullet}(G)$ such that $F_i\in A_n$ for all $i$ and
$\Delta$ is the diagonal functor then $\Delta((F_i))\in A_n$.
\item for any smooth $X$ in $G-Sch/k$ and $Z$ closed in $X$ everywhere
of codimension $\ge n$, $X/(X-Z)$ is in $A_n$
\end{enumerate}
\end{definition}
We say that an object is thick if it is 1-thick. 
\begin{remark}\rm
Note that we constructed our definition in such a way that an object
of $Spc_{\bullet}(G)$ is 0-thick if and only if it can be build out of
smooth schemes by means of homotopy colimits. In particular, unless
$k$ has resolution of singularities, it is not clear that any object
is 0-thick.
\end{remark}
\begin{lemma}
\llabel{simplprop}
The class $A_n$ of $n$-thick objects has the following
properties:
\begin{enumerate}
\item $A_n$ is closed under coproducts
\item for a morphism $f:X\sr Y$ where $X,Y$ are in $A_n$ one has
$cone(f)\in A_n$
\item for $X$ in $A_n$ one has $\Sigma^1_T X\in A_{n+1}$
\end{enumerate}
\end{lemma}
\begin{proof}
Since $A_n$ is closed under filtering colimits and any coproduct is a
filtering colimit of finite coproducts it is sufficient to show that
if $X$ and $Y$ are in $A_n$ then $X\vee Y$ is in $A_n$. Let us show
first that for any space $X$ in $A_n$, any smooth scheme $U$ and any
$Z$ of codimension at least $n+1$ in $U$, $X\vee U/(U-Z)$ is in
$A_n$. Indeed the class of $X$ for which this holds clearly satisfies
the conditions of Definition \ref{asdiag} and therefore contains
$A_n$. Similarly, the class of all $Y$ such that $X\vee Y$ is in $A_n$
clearly satisfies the first three conditions of Definition
\ref{asdiag} and we have just shown that it satisfies the
fourth. Hence it contains $A_n$. 

The cone $cone(f:X\sr Y)$ is the diagonal of a bisimplicial object
with terms $Y, Y\vee X, Y\vee X\vee X$ etc. which one obtains if one
writes the definition of the cone using second dimension for the
simplicial interval $\Delta^1_+$. Hence, any class closed under finite
coproducts and diagonals is also closed under the cones.

To verify that last condition observe that the class of $X$ such that
$\Sigma^1_T X$ is in $A_{n+1}$  satisfies the first three conditions of
Definition \ref{asdiag}. It satisfies the fourth one since
$$\Sigma^1_T (X/(X-Z))=X\times\af/(X\times\af-Z\times\{0\})$$
and codim $Z\times\{0\}=codim Z + 1$.
\end{proof}
For a finite etale $G$-scheme $W$ consider the functor $F\mapsto
F^{\wedge W}$ introduced in \cite{delnotes}.
\begin{proposition}
\llabel{nol}
Let $W$ be a finite etale $G$-scheme of degree $d$. Then the functor
$F\mapsto F^{\wedge W}$ takes $n$-thick objects to
nd-thick objects.
\end{proposition}
\begin{proof}
It is sufficient to show that the class of $F$ such that $F^{\wedge
W}$ is nd-thick satisfies the conditions of
Definition
\ref{asdiag}. The first condition follows from the fact that $F\mapsto
F^{\wedge W}$ preserves $\af$-equivalences (see \cite[p.63,
Prop. 5.2.11]{delnotes}). The second from the fact that $F\mapsto
F^{\wedge W}$ commutes with filtering colimits. The third condition is
obvious.

To see that it satisfies the fourth condition consider $X$ in
$G-Sch/k$ and $Z$ closed in $X$ everywhere of codimension $\ge n$. An
easy generalization of \cite[p.63, Rem. 5.2.8]{delnotes} shows that
$$(X/(X-Z))^{\wedge W}=X^W/(X^W-Z^W).$$
and it remains to note that if $W$ is etale then 
$codim(Z^W)=d codim(Z)$.  
\end{proof}
\begin{lemma}
\llabel{l3}
Let $X$ be a scheme with a free action of $G$ and $F$ be a pointed
solid $G$-sheaf which is n-thick. Then $Quot_G(X_+\wedge
F)$ is n-thick.
\end{lemma}
\begin{proof}
Let $R:Spc_{\bullet}(G)\sr Spc_{\bullet}(G)$ be the resolution functor
introduced in \cite[Constr. 3.6.3, p.43]{delnotes} which takes a
simplicial sheaf to a weakly equivalent one with terms being
coproducts of representable sheaves. By \cite[Prop. 5.1.4,
p.57]{delnotes} functor $Quot\circ R$ respects $\af$-equivalences and
for a solid $F$ the natural morphism $Quot(R(F))\sr Quot(F)$ is a weak
equivalence. Hence it is sufficient to show that for $F$ satisfying
the conditions of the lemma, $Quot(R(X_+\wedge F))$ is n-thick.

Since $F$ is assumed to be $n$-thick it is sufficient to
check that the class of $F$ such that $Quot(R(X_+\wedge F))$ is
$n$-thick satisfies the conditions of Definition
\ref{asdiag}. The first condition holds since the functor
$Quot(R(X_+\wedge -))$ preserves $\af$-equivalences. To see the second
one recalls that $R$ and $Quot$ commute with filtering
colimits. Similarly the third condition follows from the fact that
both $Quot$ and $R$ commute with the diagonal functor.

To check the fourth condition we have to verify that $Quot(R(X_+\wedge
(Y/(Y-Z))))$ is n-thick for any closed pair $(Y,Z)$ where
$codim(Z)\ge n$. Using again
\cite[Prop. 5.1.4, p.57]{delnotes} we see that
$$Quot(R(X_+\wedge (Y/(Y-Z))))\sr Quot(X_+\wedge (Y/(Y-Z)))$$
is an $\af$-equivalence.  We further have
$$
Quot(X_+\wedge (Y/(Y-Z)))=Quot(X\times Y)/(Quot(X\times
Y)-Quot(X\times Z))$$
Since the action of $G$ on $X$ is free, the scheme $Quot(X\times
Y)$ is smooth and since the projection $X\times Y\sr Quot(X\times
Y)$ is finite, $codim(Quot_G(X\times Z))=codim(X\times Z)=codim(Z)$.
\end{proof}

In the following two results $G$ is the trivial group.
\begin{proposition}
\llabel{stable}
Let $X$ be an $n$-thick space over a perfect field
$k$. Then the suspension spectrum $\Sigma^{\infty}_T(X)$ is in
$\Sigma^{n}_TSH^{eff}$.
\end{proposition}
\begin{lemma}
\llabel{rational}
Let $U$ be a smooth rational scheme over $k$. Then
$\tilde{\Sigma_s}(U)$ is thick.
\end{lemma}
\begin{proof}
Since $U$ is rational there is a dense open subscheme $V$ in $U$ such
that $V$ is an open subscheme in ${\bf A}^n$ for some $n$. Then
$\tilde{\Sigma_s}(V)\cong {\bf A}^n/V$ is thick and $U/V$ is
thick. We conclude by Lemma \ref{l4} and Lemma
\ref{simplprop}(2) that $\tilde{\Sigma_s}U$ is thick.
\end{proof}

\subsection{Symmetric powers of $T$}
Let $S_n$ be the symmetric group. In this section we assume that
$char(k)=0$ or $n<char(k)$. Consider the symmetric power
$$Symm^n(T^m):=Quot_{S_n}((T^m)^{\wedge n})$$
Our goal is to prove the following result.
\begin{theorem}
\llabel{main3}
For any $n\ge 2$ and $m\ge 1$ the space $Symm^n(T^m)$ is $(m+1)$-thick. 
\end{theorem}
Consider the linear representations $P(m,n)$ of $S^n$ defined by the
permutation action on $({\bf A}^m)^n$. Then
$$Symm^n(T^m)=Quot_{S_n}(Th(P(m,n)))$$
Let $V(m,n)$ be the reduced version of $P(m,n)$
$$V(m,n)=ker(p:({\bf A}^m)^n\sr {\bf A}^m)$$
where $p(x_1,\dots,x_n)=x_1+\dots+x_m$. Under our assumptions on
$char(k)$ there is an isomorphism
$$P(m,n)\cong V(m,n)\times {\bf A}^m$$
where the action on ${\bf A}^m$ is trivial. Therefore by Lemma
\ref{w1}
$$Symm^n(T^m)=Quot_{S_n}(Th(P(m,n)))=\Sigma^m_T Quot_{S_n}(Th(V(m,n)))$$
and Theorem \ref{main3} follows from Lemma \ref{simplprop}(3) and
Theorem \ref{main2} below.

Let $(n)$ be the
standard set of $n$-elements and $j_1,\dots,j_n$ integers $\ge 0$
such that 
\begin{eq}
\llabel{f6}
\sum ij_i=n.
\end{eq}
Then there is a bijection
\begin{eq}
\llabel{f5}
\phi:(n)\cong (1)^{j_1}\amalg (2)^{j_2}\amalg\dots\amalg (n)^{j_3}
\end{eq}
and any such bijection defines an embedding
$$e_{\phi}:S_1^{j_1}\times S_2^{j_2}\times\dots\times
S_n^{j_n}\sr S_n$$
The stabilizer in $S_n$ of a point $x$ in $P(m,n)$ is determined by
the set of diagonals $x_i=x_j$ which contain $x$ and subgroups of the
form $Im(e_{\phi})$ are exactly the subgroups which occur as
stabilizers of different points. 

The ajunction class of a subgroup of the form $Im(e_{\phi})$ is
determined by the integers $j_1,\dots,j_n$. The case $j_1=n$ and
$j_i=0$ for $i\ne 1$ corresponds to a point with all components being
different and the case $j_n=1$ and $j_i=0$ for $i\ne n$ to a point
with all components being the same. We choose one isomorphism of the
form (\ref{f5}) for each collection $j_1,\dots,j_n$ satisfying
(\ref{f6}) and let $H_{j_1,\dots,j_n}$ denote the corresponding
subgroup.

The same classification of stabilizers applies to $V(m,n)$ since
$V(m,n)$ is a subspace of $P(m,n)$.

The normalizer of $H=H_{j_1,\dots,j_n}$ is of the form
$G_{j_1,1}\times\dots\times G_{j_n,n}$ where $G_{j,i}$ is the
semi-direct product of $S_{j}$ and $S_{i}^{j}$ with respect to the
obvious permutational action. The quotient $N(H)/H$ is of the form
$S_{j_1}\times\dots\times S_{j_n}$.

Consider now the quotient $V(m,n)/V(m,n)^{\ge H}$. One can easily see
that it is isomorphic to the sum $\oplus V(m,i)^{\oplus j_i}$, the
action of $N(H)$ on it is the direct product of actions of $G_{j_i,i}$
on $V(m,i)^{\oplus j_i}$ and the action of $G_{j_i,i}$ is given by the
product of the standard actions of $j_i$-copies of $S_i$ on $V(m,i)$
and the permutation action of $S_{j_i}$. Summing this up we get the
following result:
\begin{lemma}
\llabel{ded}
One has an isomorphism of $N(H)/H$ spaces of the form
$$Quot_{N(H),H}(Th(V(m,n)/V(m,n)^{\ge H}))=\wedge_{i=1}^n
Quot_{S_i}(Th(V(m,i)))^{j_i}$$
\end{lemma}
\begin{theorem}
\llabel{main2}
For $n\ge 2$ and $m>0$, the space $Quot_{S_n}(Th(V(m,n)))$ is
thick.
\end{theorem}
\begin{example}\rm
For $n=2$, $V(m,n)$ is ${\bf A}^m$ with the sign action of
$S_2=Z/2$. Over a field of odd characteristic
one has
$$Quot_{S_2}(V(m,2)-\{0\})={\cal O}(-2)_{{\bf P}^{m-1}}-z({\bf
P}^{m-1})$$
where $z$ is the sero section. In particular it is a rational variety
and therefore
$$Quot_{S_2}(Th(V(m,2)))=\tilde{\Sigma_s}(Quot_{S_2}(V(m,2)-\{0\}))$$
is thick by Lemma \ref{rational}.
\end{example}
{\bf Proof of Theorem \ref{main2}:} We proceed by induction on $n\ge
2$.  Let $V_0$ be the open subscheme of $V(m,n)$ where $S_n$ acts
freely. By Lemma \ref{genpoint} below
$\tilde{\Sigma_s}Quot_{S_n}(V_0)$ is thick. In view of
Proposition \ref{basic} and Lemma \ref{simplprop}(1,2) it remains to
show that for $H=H_{j_1,\dots,j_n}$, $H\ne e,S_n$ the space
$$Quot_{N(H)}((V^{\ge H}-V^{>H})_+\wedge Quot_{N(H),H}(Th(V/V^{\ge
H})))$$
is thick. For $n=2$ any subgroup of $S_n$ is $e$ or $S_n$ and
our statement is trivial. Hence we may assume inductively that $n\ge
3$ and the theorem is proved for $n-1$. Since the action of $N(H)$ on
$V^{\ge H}-V^{>H}$ is free and the space $F=Quot_{N(H),H}(Th(V/V^{\ge
H}))$ is solid it is sufficient by Lemma \ref{l3} to see that $F$ is
thick as a $N(H)/H$-space.  Since $H\ne e, S_n$ there exists
$n>i\ge 2$ such that $j_i\ne 0$. By the inductive assumption
$Quot_{S_i}(Th(V(m,i)))$ is thick and therefore
$Quot_{S_i}(Th(V(m,i)))^{j_i}$ is thick as a $S_{j_i}$-space
by Lemma \ref{nol}. By Lemma \ref{ded}, $F$ is the smash
product of the form $F'\wedge Quot_{S_i}(Th(V(m,i)))^{j_i}$ and
$N(H)/H=G\times S_{j_i}$ where $G$ acts on the first factor and
$S_{j_i}$ on the second and we conclude that $F$ is thick.
\begin{lemma}
\llabel{genpoint}
The object $Quot_{S_n}(\tilde{\Sigma_s} V_0)$ is thick.
\end{lemma}
\begin{proof}
By Lemma \ref{rational} it is sufficient to show that
$Quot_{S_n}(V_0)$ is rational. Let $V_{0,0}$ be the open subset in
$V_0$ which consists of $x_1,\dots,x_n$ such that the first components
of all $x_i\in {\bf A}^m$ are different. The projection to the first
component 
\begin{eq}
\llabel{f7}
V_{0,0}(m,n)\sr V_0(1,n)
\end{eq}
makes $V_{0,0}(m,n)$ into an equivariant vector bundle over
$V_0(1,n)$. Since the action of $S_n$ on $V_0(1,n)$ is free, the map
of quotient schemes defined by (\ref{f7}) is a vector bundle as
well. Hence it is sufficient to show that $Quot(V_0(1,n))$ or,
equivalently, $Quot(V(1,n))$ is rational. 

The quotient $Quot(P(1,n))={\bf A}^n/S_n$ can be identified in the
standard way with ${\bf A}^n$ where the first coordinate is given on
$P(1,n)$ by the sum of components. It follows that $Quot(V(1,n))={\bf
A}^{n-1}$ and in particular that it is rational.
\end{proof}

\subsection{Reformulation for smooth schemes}
Denote for a moment the (pointed, simplicial) sheaves on smooth
schemes by $Spc_{\bullet}(Sm)$ and sheaves on all quasi-projective
schemes by $Spc_{\bullet}(Sch)$. Let $A_n$ be the class of $n$-thick
spaces in $Spc_{\bullet}(Sch)$ and $A_n^{Sm}$ the
similarly defined class in $Spc_{\bullet}(Sm)$. The functor
$$\pi_*:Spc_{\bullet}(Sch)\sr Spc_{\bullet}(Sm)$$
which takes a sheaf on $Sch/k$ to its restriction on $Sm/k$ respects
limits, colimits and $\af$-equivalences. This implies immediately that
$\pi_*$ takes $A_n$ to $A_n^{Sm}$. Therefore, all the results about
thickness proved above remain valid if we work in the context of
sheaves on smooth schemes.

\subsection{Motivic Eilenberg-MacLane spaces}

In this section $k$ is a field of characteristic zero and we work in
the context of sheaves on smooth schemes. Let $K_n=\zz_{tr}({\bf
A}^n)/\zz_{tr}({\bf A}^n-\{0\})$ be the n-th motivic Eilenberg-Maclane
space over $k$. Consider the obvious morphism $T^n\sr K_n$ where
$T^n=T^{\wedge n}={\bf A}^n/({\bf A}^n-\{0\})$.  The goal of this
section is to prove the following theorem.
\begin{theorem}
\llabel{main}
For $n>0$, the space $\Sigma_s cone(T^n\sr K_n)$ is $(n+1)$-thick.
\end{theorem}

\begin{example}\rm
Consider the case $n=1$. It is easy to see using Lemma \ref{aw} below
that $K_1$ is $\af$-equivalent to $({\bf P}^{\infty},*)$ where $*$ is
a rational point and the morphism $T\sr K_1$ corresponds to the
standard embedding $({\bf P}^1,*)\sr ({\bf P}^{\infty},*)$. The cone
is given by ${\bf P}^{\infty}/{\bf P}^1$ and its s-suspension is
2-thick by the reduced analog of Lemma \ref{l4} and Lemma
\ref{simplprop}(2).
\end{example}
Consider the sheaf $K_n^{eff}$ associated to the presheaf of the form
\begin{eq}
\llabel{presh}
K_n^{eff}:U\mapsto c^{eff}(U\times{\bf A}^n/U)/c^{eff}(U\times({\bf
A}^n-\{0\})/U)
\end{eq}
where $c^{eff}(X/U)$ is the monoid of effective finite cycles on $X$
over $U$ and the quotient on the right hand side means that two cycles
$Z_1$, $Z_2$ on $U\times{\bf A}^n$ are identified if $Z_1-Z_2$ is in
$U\times({\bf A}^n-\{0\})$. 
\begin{lemma}
\llabel{aw}
The obvious map $K^{eff}_n\sr K_n$ is an $\af$-equivalence for $n\ge
1$.
\end{lemma}
\begin{proof}
Let $F^{eff}$ be the presheaf (\ref{presh}) and let $F$ be the
presheaf given by 
$$U\mapsto \zz_{tr}({\bf A}^n)(U)/\zz_{tr}({\bf A}^n-\{0\})(U)$$
such that the sheaf associated with $F$ is $K_n$. Consider the
singular simplicial presheaves $C_*(F^{eff})$ and $C_*(F)$. The
standard argument shows that the morphisms
$$K^{eff}_n=a_{Nis}F^{eff}\sr a_{Nis}C_*(F^{eff})$$
$$K_n=a_{Nis}F\sr a_{Nis}C_*(F)$$
are $\af$-equivalences. For any smooth $U$ the abelian group
associated with the monoid $F^{eff}(U)$ coincides with
$F(U)$. Therefore, $a_{Nis}C_*(F)$ is the sheaf of simplicial abelian
groups associated with the sheaf of simplicial monoids
$a_{Nis}C_*(F^{eff})$. Let us show that the natural map 
$$a_{Nis}C_*(F^{eff})\sr a_{Nis}C_*(F)$$
is a local equivalence in the Nisnevich topology i.e. that for any
henselian local $S$ the map of simplicial sets
$$C_*(F^{eff})(S)\sr C_*(F)(S)$$
is a weak equivalence. 

For any $i$, $C_i(F^{eff})(S)=F^{eff}({\bf A}^i_S)$ is the free
commutative monoid generated by irreducible closed subsets in
$\Delta^i_S\times {\bf A}^n$ which are finite and equidimensional over
$\Delta^i_S$ and have nontrivial intersection with $\Delta^i_S\times
\{0\}$.

Therefore, by Lemma \ref{sscase} below it is sufficient to show that
$$\pi_0(C_*(F^{eff})(S))=0.$$
An element in $C_0(F^{eff})(S)$ may be represented by a cycle on ${\bf
A}^n_S$ of the form $Z=\sum n_i Z_i$ where $Z_i$ are closed
irreducible subsets of ${\bf A}^n_S$ and $n_i> 0$. Since $Z_i$ are
finite over $S$ they are local (since $S$ is henselian) and therefore
we may assumed that the closed points of all $Z_i$ lie in
$\{0\}_S$. Consider now the cycle $H$ on ${\bf A}^1_S\times{\bf A}^n$
obtained from $Z$ by the pull-back with respect to the map
$(t,x)\mapsto x-t$. This cycle is finite over $\af_S$. Its restriction
to $t=0$ is $Z$ and the restriction to $t=1$ has all the closed points
in $\{1\}_S$ and therefore lies in $({\bf A}^{n}-\{0\})_S$. Hence the
image $h$ of $H$ in $C_1(F^{eff})(S)$ has the property $\partial_0h=Z$
and $\partial_1h=0$ and we conclude that $\pi_0=0$.
\begin{lemma}
\llabel{sscase} Let $M$ be a commutative simplicial monoid and $f:M\sr
M^+$ the canonical map from $M$ to the associated simplicial group. If
for any $n$ the monoid $M_n$ is a free commutative monoid and
$\pi_0(M)$ is a group then $f$ is a weak equivalence.
\end{lemma}
\begin{proof}
If $\pi_0(K)$ is a group the topological monoid $|M|$ is a group-like
associative $H$-space. For such $H$-spaces the canonical map $|M|\sr
\Omega B(|M|)$ is a weak equivalence (follows from \cite{DoldL}). The
same argument shows that $|M^+|\sr \Omega B(|M^+|)$ is a weak
equivalence. It remains to show that $B(|M|)\sr B(|M^+|)$ is a weak
equivalence. By construction $B(|M|)$ is the same as $|B(M)|$ where
$B(M)$ is the diagonal of the bisimplicial object obtained by applying
the classifying space construction to each term of $M$. Since the class
of weak equivalences of simplicial sets is closed under taking
diagonals it remains to show that for a free commutative monoid $M$
the map $B(M)\sr B(M^+)$ is a weak equivalence. For this case see
\cite[]{}.
\end{proof}
\end{proof}
We are going now to analyze the structure of $K^{eff}_n$. Let $K^{\le
d}_n$ be the subsheaf in $K^{eff}_n$ associated with the subpresheaf
which takes $U$ to the image in $K^{eff}_n(U)$ of the set of finite
cycles on ${\bf A}^n_U/U$ of degree everywhere $\le d$ over
$U$. Clearly, $K^{eff}_n=colim_d K^{\le d}_n$. The sheaf $K^{\le 1}_n$
coincides with $T^n={\bf A}^n/({\bf A}^n-\{0\})$ naturally embedded
into $K_n^{eff}$. Consider the quotients $K^{\le d}_n/K^{\le d-1}_n$.
\begin{lemma}
\llabel{quot}
Let $k$ be a field of characteristic zero. Then there is an
isomorphism 
\begin{eq}
\llabel{f8}
K^{\le d}_n/K^{\le d-1}_n\cong Symm^d(T^n)
\end{eq}
where $Symm^d(T^n):=Quot_{S_d}((T^n)^{\wedge d})$.
\end{lemma}
\begin{proof}
Let $c^{\le d}({\bf A}^n)$ be the sheaf of finite cycles of degree
$\le d$ on ${\bf A}^n$. By \cite[]{}, for $k$ of characteristic zero
this sheaf is represented by 
$$\coprod_{i\le d} Symm^i({\bf A}^n)=\coprod_{i\le d} Quot_{S_i}({\bf
A}^{ni})$$
Since the right hand side of (\ref{f8}) is given by $Quot_{S_d}({\bf
A}^{nd}/({\bf A}^{nd}-\{0\}))$ it remains to show that the restriction
of the map 
$$c^{\le d}({\bf A}^n)\sr K^{\le d}_n/K^{\le d-1}_n$$
to $Quot_{S_d}({\bf A}^{nd})$ is surjective and that two sections
of $Quot_{S_d}({\bf A}^{nd})$ over a henselian local scheme $S$ project to
the same section of $K^{\le d}_n/K^{\le d-1}_n$ if and only if they
project to the same section of $Quot_{S_d}({\bf A}^{nd}/({\bf
A}^{nd}-\{0\}))$.

The surjectivity is clear. Let $Z=\sum n_i Z_i$, $W=\sum m_j W_j$ be
two sections of $c^{\le d}({\bf A}^n)$ over $S$ of degree strictly
$d$. By definition they project to the same section of $K^{\le
d}_n/K^{\le d-1}_n$ if both $Z$ and $W$ are equivalent modulo the
components lying in ${\bf A}^n-\{0\}$ to cycle of degree $\le
d-1$. Since $Z_i$ and $W_j$ are local this means that not all of the
closed points of the $Z_i$'s (resp. $W_j$'s) lie in $\{0\}$. This is
equivalent to the condition that $Z$ and $W$ lie in $Quot_{S_d}({\bf
A}^{nd}-\{0\})$.
\end{proof}
{\bf Proof of Theorem \ref{main}:} By Lemma \ref{aw}, it is sufficient
to show that the space ${\Sigma_s}(K^{eff}_n/T^n)$ is
$(n+1)$-thick. Since $T^n=K^{\le 1}_n$ the space $K^{eff}_n/T^n$ has a
filtration which starts with $K^{\le 2}_n/K^{\le 1}_n$ and has
quotients of the form $K^{\le d}_n/K^{\le d-1}_n$ for $d\ge 3$. Since
the class of thick objects is closed under filtering colimits
it is sufficient by the obvious reduced analog of Lemma \ref{l4} and
Lemma \ref{simplprop}(2) to show that $K^{\le d}_n/K^{\le d-1}_n$ are
$(n+1)$-thick  for $d\ge 2$. This follows immediately
from Lemma
\ref{quot} and Theorem \ref{main3}.

As a corollary of Theorem \ref{main} we can prove the main conjecture
of \cite{Open} over fildes of characteristic zero.
\begin{theorem}
\llabel{appl}
Let $k$ be a field of characteristic zero. Then $s_0({\bf
1})=H_{\zz}$.
\end{theorem}
\begin{proof}
Consider the unit map $e:{\bf 1}\sr H_{\zz}$. Since
$s_0(H_{\zz})=H_{\zz}$ it is sufficient to show that $s_0(e)$ is an
isomorphism i.e. that $s_0(cone(e))=0$. By definition of $s_0$ this
would follow if we can show that $cone(e)$ belongs to
$\Sigma^{1}_TSH^{eff}$. The map $e$ is given on the level of spectra
by the maps $T^n\sr K_n$ considered above. In particular we have
$$cone(e)=hocolim_n \Sigma^{-n}_T\Sigma^{\infty}_T(cone(T^n\sr K_n))$$
Since $\Sigma^{1}_TSH^{eff}$ is closed under homotopy colimits it
remains to check that
$$\Sigma^{-n}_T\Sigma^{\infty}_T(cone(T^n\sr
K_n))\in\Sigma^{1}_TSH^{eff}$$
i.e. that
$$\Sigma^{\infty}_T(cone(T^n\sr
K_n))\in  \Sigma^{n+1}_TSH^{eff}$$
Since $\Sigma^{n+1}_TSH^{eff}$ is stable under $\Sigma_s$ it follows
from Theorem \ref{main} and Proposition \ref{stable}.
\end{proof}

\end{document}